\documentclass[makeidx]{amsart}

\usepackage{a4wide}
\usepackage{amssymb}
\usepackage[latin1]{inputenc}
\usepackage{graphicx}
\usepackage{dsfont}
\usepackage{color}
\usepackage[hyperref]{hyperref}

\renewcommand{\phi}{\varphi}
\newcommand{\C}{{\mathbb{C}}}

\newcommand{\R}{{\mathbb{R}}}

\newcommand{\Z}{{\mathbb{Z}}}

\renewcommand{\epsilon}{\varepsilon}
\renewcommand{\theta}{\vartheta}
\renewcommand{\S}{{\mathbb{S}}}
\newcommand{\T}{{\mathbb{T}}}
\newcommand{\Disk}[1]{{\mathbb{D}^{#1}}}

\newcommand{\pairing}[2]{{\langle{#1}|{#2}\rangle}}
\newcommand{\plastikstufe}[1]{{\mathcal{PS}({#1})}}
\newcommand{\norm}[1]{{\lVert #1\rVert}}
\newcommand{\abs}[1]{{\left\lvert #1\right\rvert}}

\newcommand{\overtwisted}{{\mathbb{D}_\mathrm{OT}}}
\newcommand{\standardalpha}{{\alpha_{\mathrm{std}}}}
\newcommand{\z}{{\mathbf{z}}}
\newcommand{\lcan}{{\lambda_{\mathrm{can}}}}

\DeclareMathOperator{\glue}{Glue}

\DeclareMathOperator{\SO}{SO}

\theoremstyle{plain}
\newtheorem{theorem}{Theorem}
\newtheorem{propo}[theorem]{Proposition}
\newtheorem{coro}[theorem]{Corollary}

\theoremstyle{remark}
\newtheorem{remark}[theorem]{Remark}

\theoremstyle{definition}
\newtheorem*{defi}{Definition}

\numberwithin{equation}{section}

\begin{document}

\title{Every contact manifold can be given a non-fillable contact
  structure}

\author[K.~Niederkrüger]{Klaus Niederkrüger}
\author[O.~van~Koert]{Otto van Koert}

\email[K.~Niederkrüger]{kniederk@ulb.ac.be}
\email[O.~van~Koert]{ovkoert@ulb.ac.be}

\address{Département de Mathématiques, Université Libre de Bruxelles, CP 218\\
  Boulevard du Triomphe\\B-1050 Bruxelles\\Belgium}

\begin{abstract}
  Recently Francisco Presas~Mata constructed the first examples of
  closed contact manifolds of dimension larger than $3$ that contain a
  plastikstufe, and hence are non-fillable.  Using contact surgery on
  his examples we create on every sphere $\S^{2n-1}$, $n\geq 2$, an
  exotic contact structure $\xi_-$ that also contains a plastikstufe.  As a
  consequence, every closed contact manifold
  $(M,\xi)$ (except $\S^1$) can be converted into a contact
  manifold that is not (semi-positively) fillable by taking the
  connected sum $(M,\xi)\# (\S^{2n-1},\xi_-)$.
\end{abstract}


\maketitle

Most of the natural examples of contact manifolds can be realized as
convex boundaries of symplectic manifolds.  These manifolds are called \emph{symplectically fillable}. An important class of
contact manifolds that do not fall into this category are so-called
\emph{overtwisted} manifolds (\cite{Eliashberg_NonFillable},
\cite{Gromov_Kurven}).  Unfortunately, the notion of overtwistedness
is only defined for $3$--manifolds.  A manifold is overtwisted if one
finds an embedded disk $\Disk{2}$ such that
$\left.T\Disk{2}\right|_{\partial \Disk{2}}\subset \xi$, an
\emph{overtwisted} disk $\overtwisted$. This topological definition gives an effective way to
find many examples of contact $3$--manifolds that are non-fillable.

Until recently no example of a non-fillable contact manifold
in higher dimension was known, but Francisco Presas Mata recently
discovered a construction that allowed him to build many non-fillable
contact manifolds of arbitrary dimension
(\cite{PresasExamplesPlastikstufes}).  He showed that after performing
this construction on certain manifolds, they admit the embedding of a
\emph{plastikstufe}.  Roughly speaking, a plastikstufe
$\plastikstufe{S}$ can be thought of as a disk-bundle over a closed
$(n-2)$--dimensional submanifold $S$, where each fiber looks like an
overtwisted disk.  As shown in~\cite{NiederkruegerPlastikstufe}, the
existence of such an object in a contact manifold $(M^{2n-1},\xi)$
excludes the existence of a symplectic filling.

In this paper we extend Presas' results to a much larger class of
contact manifolds. The idea is to start with one of his examples and
to use contact surgery to simplify the topology and convert it into a
contact sphere\footnote{In the scope of this article a \emph{contact
    sphere} will be a smooth sphere carrying a contact structure.}
$S_E$ that admits the embedding a plastikstufe. If $(M,\xi)$ is any
other contact manifold, then $(M,\xi)\# S_E \cong M$ carries a contact
structure that also has an embedded plastikstufe and is hence
non-fillable.

\begin{defi}
  Let $(M,\alpha)$ be a cooriented $(2n-1)$--dimensional contact
  manifold, and let $S$ be a closed $(n-2)$--dimensional manifold.  A
  \textbf{plastikstufe $\plastikstufe{S}$ with singular set $S$} in
  $M$ is an embedding of the $n$--dimensional manifold
  \begin{equation*}
    \iota:\,\Disk{2}\times S \hookrightarrow M
  \end{equation*}
  that carries a (singular) Legendrian foliation given by the
  $1$--form $\beta := \iota^*\alpha$ satisfying:
  \begin{itemize}
  \item The boundary $\partial\plastikstufe{S}$ of the plastikstufe is
    the only closed leaf.
  \item There is an elliptic singular set at $\{0\}\times S$.
  \item The rest of the plastikstufe is foliated by an
    $\S^1$--family of stripes, each one diffeomorphic to $(0,1)\times
    S$, which are spanned between the singular set on one end and
    approach $\partial\plastikstufe{S}$ on the other side
    asymptotically.
  \end{itemize}
\end{defi}

The importance of the plastikstufe lies in the following theorem.

\begin{theorem}
  Let $(M,\alpha)$ be a contact manifold containing an embedded
  plastikstufe.  Then $M$ does not have a semipositive strong symplectic
  filling.  In particular, if $\dim M \le 5$, then $M$ does not have any strong symplectic
  filling at all.
\end{theorem}

\begin{defi}
  A contact manifold is called \textbf{$PS$--overtwisted} if it admits
  the embedding of a plastikstufe.
\end{defi}

Whether $PS$--overtwisted can be taken as the general definition of
overtwisted in higher dimensions, has yet to be clarified in the
future (see Remark~\ref{bemerkung ueber ueberdreht versionen}).

In dimension~$3$, the definition of $PS$--overtwisted is identical to
the standard definition of overtwistedness.  Using the Lutz twist, it
is easy to convert a tight contact structure on a $3$--manifold $M$
into an overtwisted one.  Until recently, no example of a closed
$PS$--overtwisted contact manifold of dimension larger than~$3$ was
known, but Francisco Presas Mata found a beautiful construction, which
allowed him to create such examples in arbitrary
dimension~\cite{PresasExamplesPlastikstufes}.

\begin{theorem}[F.~Presas Mata]\label{satz: presas gluing} Let
  $(M,\alpha)$ be a contact manifold, which contains a
  $PS$--overtwisted contact submanifold $N$ of codimension~$2$.
  Assume that $N$ has trivial normal bundle.  Then we can glue
  $N\times\T^2$ via a fiber sum to $M$, and the resulting manifold
  \begin{equation*}
    \glue(M,N,\alpha) := M \cup_N N\times\T^2
  \end{equation*}
  supports a $PS$--overtwisted contact structure that coincides on $M$
  outside a small neighborhood of the gluing area with the original
  structure.  More precisely, if $N$ contains a plastikstufe
  $\plastikstufe{S}$, then $\glue(M,N,\alpha)$ contains
  $\plastikstufe{S\times\S^1}$.
\end{theorem}

With this construction, he was immediately able to find
$PS$--overtwisted contact manifolds in every odd dimension greater than $1$.

\begin{coro}\label{beispiele von presas}
  \begin{itemize}
  \item [(1)] There is a contact form $\alpha_-$ on the $5$--sphere
    $\S^5$ that is obtained by taking an open book with a single
    left-handed Dehn-twist (see Section~\ref{sec:link-haendige
      sphaeren}). It restricts on the standard embedding of $\S^3$ to
    an overtwisted contact structure.  Hence, one finds on the
    manifold
    \begin{equation*}
      M_0 = \glue(\S^5,\S^3,\alpha_-) = \S^5 \cup_{\S^3} \S^3\times \T^2
    \end{equation*}
    a $PS$--overtwisted contact structure.
  \item [(2)] Let $\Sigma_g$ denote the closed Riemann surface of
    genus~$g\ge 2$, and let $(M,\alpha)$ be a closed $PS$--overtwisted
    contact manifold.  Then the manifold $M\times \Sigma_g$ also
    supports a $PS$--overtwisted contact structure.
  \end{itemize}
\end{coro}

In this note, we apply contact surgery to examples that are similar to
the manifold $M_0$, and we obtain the following corollary.

\begin{coro}\label{hauptfolgerung}
  Every sphere $\S^{2n+1}$ with $n\geq 1$ supports a
  $PS$--overtwisted contact structure.  More precisely, on
  $\S^{2n+1}$, with $n\ge 1$, exists a contact structure which admits
  the embedding of a plastikstufe $\plastikstufe{\T^{n-1}}$ (with
  $\T^0 := \{p\}$ and $\T^1 := \S^1$).
\end{coro}

\begin{coro}
  It is possible to modify the contact structure $\xi = \ker\alpha$ of
  a contact manifold $(M,\alpha)$ with $\dim M\ge 3$ in an arbitrary
  small open set in such a way that the new contact structure is
  $PS$--overtwisted.
\end{coro}
\begin{proof}
  Attach one of the contact spheres obtained in
  Corollary~\ref{hauptfolgerung} via connected contact sum to the
  manifold~$M$.
\end{proof}

\begin{remark}\label{bemerkung ueber ueberdreht versionen} The results
  above make it tempting to claim that $PS$--overtwisted is the proper
  definition of overtwistedness in higher dimensions.  The second
  author proved together with Frédéric Bourgeois that contact
  manifolds having an open book decomposition of a certain type have
  vanishing contact homology, which as \emph{folklore} tells also
  implies that they are non-fillable.  In $3$~dimensions, overtwisted
  can equivalently be defined via open book decompositions, and thus
  in higher dimensions it would also be interesting to find the
  precise relation between the definition using a plastikstufe or an
  open book decomposition.

  In particular all of the spheres $(\S^{2n-1},\alpha_-)$ defined in
  Section~\ref{sec:link-haendige sphaeren} have vanishing contact
  homology, and hence are non-fillable even before applying the Presas
  gluing.
\end{remark}

\subsection*{Acknowledgments}

This article was written at the \emph{Université Libre de Bruxelles},
where we are both being funded by the \emph{Fonds National de la
  Recherche Scientifique} (FNRS).  We thank Hansjörg Geiges for
fruitful discussions.

\section{The contact spheres
  $(\S^{2n-1},\alpha_-)$}\label{sec:link-haendige
sphaeren}

In this section, we will describe an exotic contact structure on each
sphere $\S^{2n-1}$ with $n\geq 2$, with the interesting property that they can
all be stacked into each other in a natural way.

Let $\S^{2n-1}$ be the unit sphere in~$\C^n$ with coordinates $\z =
(z_1,\dotsc,z_n)\in\C^n$, and let $f$ be the polynomial
\begin{equation*}
  f:\,\C^n\to\C,\, (z_1,\dotsc,z_n) \mapsto z_1^2+ \dotsb +z_n^2 \;.
\end{equation*}
The $1$--form
\begin{align*}
  \alpha_- & := i\, \sum_{j=1}^n\left(z_j\,d\bar z_j - \bar
    z_j\,dz_j\right) - i\, \left( f\,d \overline{f} - \overline{f}\,d
    f \right)
\end{align*}
defines a contact structure on $\S^{2n-1}$, which is obtained by using
an open book decomposition with left-handed Dehn-twist (see
Remark~\ref{remark: open book}).  The first term of $\alpha_-$ is just
the standard contact form $\standardalpha$ on the sphere, and it is
the second term that is responsible for changing the properties
of~$\alpha_-$.

\begin{propo}
  The sphere $(\S^{2n-1},\alpha_-)$ is a contact manifold.
\end{propo}
\begin{proof}
  Consider the $1$--form
  \begin{align*}
    \widetilde\alpha_- & := i\, \sum_{j=1}^n\left(z_j\,d\bar z_j -
      \bar z_j\,dz_j\right) - i\, \left( \frac{f}{\norm{\z}}\,d
      \frac{\overline{f}}{\norm{\z}} -
      \frac{\overline{f}}{\norm{\z}}\,d \frac{f}{\norm{\z}} \right)\;.
  \end{align*}
  Its restriction to the unit sphere is equal to $\alpha_-$.  We will
  show that its exterior differential, the $2$--form
  \begin{align*}
    \omega_- & := d\widetilde\alpha_- = 2i\, \sum_{j=1}^n d z_j\wedge
    d\bar z_j - 2i \, d\;\left(\frac{f}{\norm{\z}}\right)\wedge
    d\;\left(\frac{\overline{f}}{\norm{\z}}\right) \;,
  \end{align*}
  is symplectic on $\C^n-\{0\}$.  To compute the $n$--fold product
  $\omega_-^n$, note that the last term can appear at most once in
  each term of the total product.  Furthermore since the first terms
  always couple a~$dz_j$ with $d\bar z_j$, this eliminates most mixed
  forms of the last term.  Using this, one easily computes
  \begin{align*}
    \omega_-^n &= -\frac{(2i)^n n!}{\norm{\z}^4}\,
    \left(3\,\norm{\z}^{4} -2\,\abs{f}^2\right)\,dz_1\wedge d\bar
    z_1\wedge \dotsb \wedge dz_n\wedge d\bar z_n\;.
  \end{align*}
  This form does not vanish, because
  \begin{equation*}
    \abs{f}^2 =  \abs{z_1^2+ \dotsb +z_n^2}^2 =
    \abs{\pairing{\bar\z}{\z}}^2 
    \le \norm{\z}^4\;.
  \end{equation*}
  Now, it follows that $\widetilde\alpha_-$ (and hence also
  $\alpha_-$) is a contact form by using that the Liouville field
  \begin{equation*}
    X_L = \frac{1}{2}\,\partial_r = \frac{1}{2}\,(z_1,\dotsc,z_n)\;.
  \end{equation*}
  for $\omega_-$ satisfies $\iota_{X_L}\omega_- = \widetilde\alpha_-$.
\end{proof}

\begin{remark}
\label{rem:overtwistedforms}
  \begin{enumerate}
  \item We obtain sequences of contact embeddings
    \begin{equation*}
      (\S^3,\alpha_-) \hookrightarrow (\S^5,\alpha_-) \hookrightarrow
      (\S^7,\alpha_-) \hookrightarrow \dotsm \;,
    \end{equation*}
    where each map is just given by 
\begin{align*}
  \begin{aligned}
    \iota_j:\, \S^{2k-1} &\longrightarrow \S^{2k+1},\, \\
    (z_1,\dotsc,z_k) & \longmapsto
    (z_1,\dotsc,z_{j-1},0,z_j,\dotsc,z_k) \;.
  \end{aligned}
\end{align*}

    The normal bundle of every $(2k-1)$--sphere in the following
    $(2k+1)$--sphere is of course trivial.
  \item \label{remark: open book} Using similar computations as the
    ones in \cite{KoertVan} or \cite{NiederkruegerSO3Contact}, it can
    be seen that $(\S^{2n-1}, \alpha_-)$ is compatible with the open
    book $(B = f^{-1}(0), \theta = \overline{f}/\abs{f})$, which is
    equivalent to the abstract open book with page
    $(T^*\S^{n-1},d\lcan)$ and monodromy map consisting of a single
    left-handed Dehn-twist.  In particular, the $3$--dimensional case
    $(\S^3,\alpha_-)$ is overtwisted, and it is not difficult to
    localize an overtwisted disk: The intersection $F$ between the
    $3$--sphere $\S^3$ and the hyperplane
    $\bigl\{(x_1+iy_1,x_2+iy_2)\bigm| y_1=x_2\bigr\}$ is diffeomorphic
    to a $2$--sphere, which is foliated by $\alpha_-$.  Using
    stereographic projection
    \begin{align*}
      \Phi:\,& \C \to \C^2 \\
      & x+iy \mapsto \frac{1}{\sqrt{2}\,(1+x^2+y^2)} \,
      \begin{pmatrix}
        (x+1)^2 + y^2 -2 \\ 2y \\ 2y \\ (x-1)^2 + y^2 - 2
      \end{pmatrix}\;,
    \end{align*}
    we obtain for the pull-back of the contact form onto this sphere
    \begin{align*}
      \Phi^*\alpha_- & = \frac{4\,(3r^4-10r^2+3\bigr)}{(1+r^2)^4}\,\bigl(
      y\,dx - x\,dy\bigr)\;,
    \end{align*}
    with $r^2 = \abs{z}^2 = x^2 + y^2$.  This form does not vanish
    with exception of the points corresponding to the origin and the
    circles of radius $1/\sqrt{3}$ and $\sqrt{3}$.  Hence we find an
    overtwisted disk in each of the hemispheres of $F$ (the set
    $\Disk{2}_{\abs{z}^2 \le 1/3} := \{z\in\C |\; \abs{z}^2 \le 1/3\}$
    and $\{z\in\C |\; \abs{z}^2 \ge 3\}\cup \{\infty\}$).

  \end{enumerate}
\end{remark}

\section{$\S^5$ supports a $PS$--overtwisted contact
  structure}\label{sec: beweis dimension 5}

In this section, we will prove the result stated in
Corollary~\ref{hauptfolgerung} for dimension~$5$.  To achieve our
goal, we will simply start with the manifold $M_0$ (see
Corollary~\ref{beispiele von presas}), and then use contact
$1$--surgery to kill the fundamental group.  The general proof in
Section~\ref{sec:hauptbeweis} includes the $5$--dimensional case, but
the induction used there is relatively complicated, so that we
preferred to work out this case explicitly.  The following
proposition is well known to topologists, but since it is key
to our construction we include a proof.

\begin{propo}\label{mache einfach zusammenhaengend} Let $M$ be an
  orientable manifold of dimension~$n >3$.  Assume the fundamental
  group~$\pi_1(M)$ to be generated by the closed embedded paths
  $\gamma_1, \dotsc, \gamma_N$.  Then by using surgery on these
  circles we obtain a simply connected manifold $\widetilde M$.
\end{propo}
\begin{proof}
  We have to show that $\pi_1(\widetilde M)$ vanishes.  First note
  that the statement in the proposition needs a little clarification.
  We want all generators $\gamma_1,\dotsc,\gamma_N$ to be disjoint,
  but this is strictly speaking not possible, because all elements in
  $\pi_1(M,p_0)$ have at least the base point~$p_0\in M$ in common.
  Instead move each of the circle with a small isotopy to make them
  disjoint from each other.  Choose now new representatives
  $\gamma_1^\prime, \dotsc, \gamma_N^\prime \in \pi_1(M, p_0)$ such
  that $\gamma_j^\prime$ consists of a short segment connecting $p_0$
  with a point $(\gamma_j(0),x)$ on the boundary of the tubular
  neighborhood $\S^1\times\Disk{n-1}$ of $\gamma_j$, the path
  $(\gamma_j,x) \subset \partial(\S^1\times \Disk{n-1})$, and a copy
  of the first segment but with opposite orientation connecting
  $(\gamma_j(0),x) = (\gamma_j(1),x)$ back with $p_0$.  After the
  surgery, each of the $\gamma_j^\prime$ can be contracted to the
  point~$p_0$.

  \begin{figure}[htbp]
    \begin{picture}(0,0)%
      \includegraphics{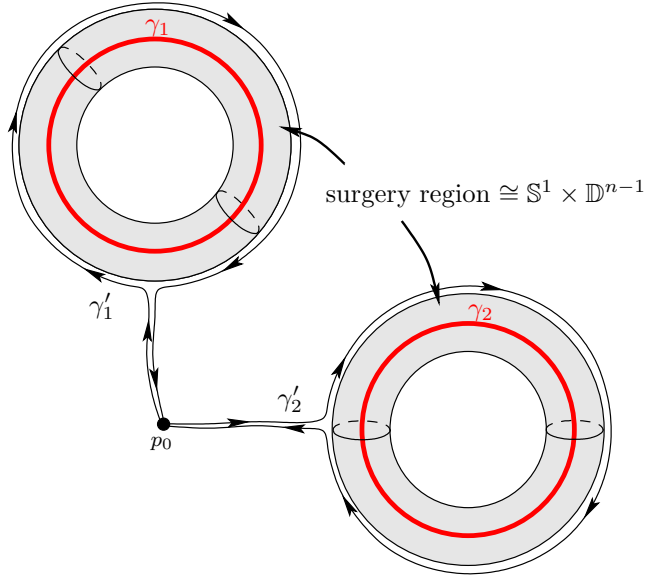}%
    \end{picture}%
    \setlength{\unitlength}{3947sp}%
    \begingroup\makeatletter\ifx\SetFigFont\undefined%
    \gdef\SetFigFont#1#2#3#4#5{%
      \reset@font\fontsize{#1}{#2pt}%
      \fontfamily{#3}\fontseries{#4}\fontshape{#5}%
      \selectfont}%
    \fi\endgroup%
    \begin{picture}(3806,3613)(1676,-5426)
      \put(3678,-3076){\makebox(0,0)[lb]{\smash{{\SetFigFont{10}{12.0}{\familydefault}{\mddefault}{\updefault}{\color[rgb]{0,0,0}surgery
                region $\cong \S^1\times \Disk{n-1}$}%
            }}}}

      \put(2579,-4612){\makebox(0,0)[lb]{\smash{{\SetFigFont{8}{9.6}{\familydefault}{\mddefault}{\updefault}{\color[rgb]{0,0,0}$p_0$}%
            }}}}

      \put(4572,-3792){\makebox(0,0)[lb]{\smash{{\SetFigFont{10}{12.0}{\familydefault}{\mddefault}{\updefault}{\color[rgb]{1,0,0}$\gamma_2$}%
            }}}}

      \put(2545,-2003){\makebox(0,0)[lb]{\smash{{\SetFigFont{10}{12.0}{\familydefault}{\mddefault}{\updefault}{\color[rgb]{1,0,0}$\gamma_1$}%
            }}}}

      \put(3375,-4355){\makebox(0,0)[lb]{\smash{{\SetFigFont{10}{12.0}{\familydefault}{\mddefault}{\updefault}{\color[rgb]{0,0,0}$\gamma_2^\prime$}%
            }}}}

      \put(2193,-3771){\makebox(0,0)[lb]{\smash{{\SetFigFont{10}{12.0}{\familydefault}{\mddefault}{\updefault}{\color[rgb]{0,0,0}$\gamma_1^\prime$}%
            }}}}
    \end{picture}%
    \caption{Every loop in the surgered space $\widetilde M$ is
      homotopic to a concatenation of circles
      $\gamma_1,\dotsc,\gamma_N$.  The surgery makes each of these
      circles contractible so that $\widetilde M$ is simply
      connected.}
    \label{bild: chirurgie fuer reduktion der fundamentalgruppe}
  \end{figure}

  Let $\gamma$ be a closed path that represents an element in
  $\pi_1(\widetilde M,p_0)$, where we assume that the base point $p_0$
  lies outside the surgery area.  With a homotopy, we can make it also
  everywhere disjoint from the surgery area which is essentially an
  $(n-2)$--sphere.  This way we obtain a loop that lives not only in
  $\widetilde M$ but also in $M$, and represents an element in the
  fundamental group $\pi_1(M, p_0)$.  In~$M$, this circle is homotopic
  to a product of the $\gamma_1^\prime,\dotsc, \gamma_N^\prime$, and
  this homotopy can be made disjoint from the surgery regions
  $\S^1\times \Disk{n-1}$, because a homotopy of curves is a map
  $[0,1]\times\S^1 \to M$, but in an $n$--dimensional manifold with
  $n>3$, it is always possible to make a the image of a $2$--manifold
  by a perturbation disjoint from a $1$--dimensional submanifold.
  This homotopy can thus be also realized in $\widetilde M$, and so
  $\gamma$ is also in $\widetilde M$ homotopic to a product of the
  $\gamma_1^\prime, \dotsc, \gamma_N^\prime$, which are all
  contractible in~$\widetilde M$.  If follows that $\gamma$ represents
  the trivial element in $\pi_1(\widetilde M)$, and hence
  $\pi_1(\widetilde M) = \{0\}$.
\end{proof}

We start by using the manifold $M_0$ given in Corollary~\ref{beispiele
  von presas}, part~(1).  The homology\footnote{From now on, we will always
  assume \emph{integer coefficients} for the homology groups.} of this
space is (as already stated by Presas) $H_0(M_0) \cong H_5(M_0) \cong
\Z$, $H_1(M_0) \cong H_4(M_0) \cong \Z^2$ and $H_2(M_0) \cong H_3(M_0)
= \{0\}$.  The fundamental group $\pi_1(M_0)$ can be easily computed
using the Seifert-van~Kampen theorem.  We obtain
\begin{equation*}
  \pi_1(M_0) = \langle a,b,c|\, aba^{-1}b^{-1} = c \rangle =
  \langle a,b \rangle = \Z * \Z \;,
\end{equation*}
where $a,b$ are the generators of $\pi_1(\S^3\times (\T^2-\Disk{2}))
\cong \pi_1(\T^2 - \Disk{2}) \cong \Z*\Z$, and $c$ generates the
fundamental group of $\S^5-\S^3\times\Disk{2} \cong
\Disk{4}\times\S^1$.  The relation follows from identifying elements
in the intersection $\S^3\times \S^1$.

Represent the two generators $a,b$ by smooth embedded paths $\gamma_1,
\gamma_2\subset \S^3\times(\T^2-\Disk{2})$ of the form
\begin{align*}
  \gamma_1:\,& [0,1] \to \S^3\times(\T^2 - \Disk{2}), \, t\mapsto
  (p_1;e^{2\pi i t},1) \\
  \gamma_2:\,& [0,1] \to \S^3\times(\T^2 - \Disk{2}), \, t\mapsto
  (p_2;1,e^{2\pi i t})
\end{align*}
with fixed $p_1\ne p_2\in\S^3$ such that both curves are isotropic and
do not intersect the plastikstufe $\plastikstufe{\S^1}$ lying in
$M_0$.  This can be achieved by choosing both points $p_1,p_2$ in the
$1$--dimensional binding of $(\S^3,\alpha_-)$ (because then both
functions $f_1,f_2$ used for the definition of $\alpha_\epsilon$ on
$\S^3\times\T^2$ vanish, cf.\ Section~\ref{sec: presas gluing auf
  untermannigfaltigkeiten}).  One can find an overtwisted disk in
$\S^3$ that intersects the binding at only one point.  For the
construction of $\plastikstufe{\S^1}$, this disk is transported
through $\S^3\times\T^2$ along two paths which are parallel to
$\gamma_1$ or $\gamma_2$.  Hence there is sufficient space to choose
$p_1,p_2$ in such a way that $\gamma_1$ and $\gamma_2$ do not
intersect $\plastikstufe{\S^1}$.

Now apply contact surgery on these generators of $\pi_1(M_0)$, i.e.\
cut out a tubular neighborhood of the two isotropic curves
representing $a$ and $b$.  This neighborhood is of the form
$\S^1\times \Disk{4} \dot\cup \S^1\times\Disk{4}$ and it has boundary
$\S^1\times\S^{3} \dot\cup \S^1\times\S^{3}$.  Glue in two copies of
$\Disk{2}\times \S^{3}$, which have the same boundary as the cavities.
As explained in \cite{WeinsteinSurgery}, the manifold $\widetilde M_0$
we obtain this way, carries a contact structure, which coincides
outside the surgery loci with the original structure, so in particular
$\widetilde M_0$ still contains a plastikstufe.  By
Proposition~\ref{mache einfach zusammenhaengend}, it follows that
$\widetilde M_0$ is simply connected.

\subsection{The homology of $\widetilde M_0$}

We want to show that $H_*(\widetilde M_0) \cong H_*(\S^5)$, because
this together with $\pi_1(\widetilde M_0) = \{0\}$, using the Poincaré
conjecture proved by Smale and the non-existence of exotic $5$--spheres shows that
$\widetilde M_0$ is diffeomorphic to $\S^5$.

All of the computations in this section are standard applications of
the Mayer-Vietoris sequence.  Use the following notation $B =
\S^1\times \Disk{4}\dot\cup\S^1\times \Disk{4}$, $\widetilde B =
\Disk{2}\times \S^3 \dot\cup \Disk{2}\times \S^3$, $A = M_0 - B =
\widetilde M_0 - \widetilde B$, $A\cap B = A\cap \widetilde B =
\S^1\times \S^3 \dot\cup \S^1\times \S^3$.  Then, because $H_2(A\cap
B) = \{0\}$, both Mayer-Vietoris sequences for $M_0 = A\cup B$ and for
$\widetilde M_0 = A\cup \widetilde B$ split at that homology group.
Using that $H_2(B) \cong H_2(\widetilde B) \cong H_1(\widetilde M_0) =
\{0\}$, the Mayer-Vietoris sequences for the pairs $(A,B)$ and
$(A,\widetilde B)$ reduce to
\begin{align*}
  0 \rightarrow& H_2(A) \oplus H_2(B) \rightarrow H_2(M_0) \rightarrow
  H_1(A\cap B) \rightarrow H_1(A)\oplus H_1(B) \rightarrow H_1(M_0)
  \rightarrow 0 \\
  0 \rightarrow& H_2(A) \oplus H_2(\widetilde B) \rightarrow
  H_2(\widetilde M_0) \rightarrow H_1(A\cap \widetilde B) \rightarrow
  H_1(A)\oplus H_1(\widetilde B) \rightarrow 0 \;.
\end{align*}
Because $H_2(M_0) = \{0\}$, it follows that $H_2(A)$ also vanishes.
The top sequence simplifies to
\begin{align*}
  0 \rightarrow \Z^2 \rightarrow H_1(A)\oplus \Z^2 \rightarrow \Z^2
  \rightarrow 0 \;.
\end{align*}
As $H_1(A)$ cannot have torsion in such a short exact sequence, we
obtain $H_1(A)\cong \Z^2$.  With this, the second sequence simplifies
to
\begin{align*}
  0 \rightarrow H_2(\widetilde M_0) \rightarrow \Z^2 \rightarrow \Z^2
  \rightarrow 0 \;,
\end{align*}
so that $H_2(\widetilde M_0)$ vanishes.

By Poincaré duality and the universal coefficient theorem, we have
$H_3(\widetilde M_0) \cong H_4(\widetilde M_0) = \{0\}$. This implies
that $\widetilde M_0$ is homeomorphic to $\S^5$ and in fact it is even
diffeomorphic to $\S^5$, because there are no exotic $5$--spheres.
This proves Corollary~\ref{hauptfolgerung} for dimension~$5$.

\section{Submanifolds and Presas gluing}
\label{sec: presas gluing auf untermannigfaltigkeiten}

In \cite{PresasExamplesPlastikstufes}, Presas starts out with a
manifold $(M,\alpha)$ which contains a codimension~$2$ contact
submanifold $N$ with trivial normal bundle.  A neighborhood of $N$ in
$M$ is contactomorphic to $(N\times \Disk{2}, \left.\alpha\right|_{TN}
+ r^2\,d\phi)$.  The product manifold $N\times\T^2$ also supports a
contact structure, namely $\alpha_\epsilon := \left.\alpha\right|_{TN}
+ \epsilon\,\left(f_1\,d\theta_1 + f_2\,d\theta_2\right)$, where
$(e^{i\theta_1},e^{i\theta_2})$ are the coordinates of the $2$--torus
and $f_1,f_2:\,N\to \R$ are certain functions, which are obtained from
an open book decomposition of $N$ (see \cite{BourgeoisTori}).  A fiber
$N\times\{p_0\}\subset N\times \T^2$ has a neighborhood that is
contactomorphic to $(N\times\Disk{2}, \left.\alpha\right|_{TN} -
r^2\,d\phi)$, and one can perform the fiber connected sum of
$N\times\T^2$ onto $M$ along $N\times\{p_0\}$ and $N\subset M$.
Presas has shown that if $N$ is $PS$--overtwisted then so is $M\cup_N
N\times\T^2$.  This construction can be carried out simultaneously on
several embedded contact manifolds.  To this end we need a
neighborhood theorem that is adapted to such a situation.

\begin{propo}\label{umgebungssatz fuer zwei mgktn} Let $N, S_1,\dotsc,
  S_r$ be codimension~$2$ contact submanifolds of $(M,\alpha)$ with
  trivial normal bundle.  Assume that all of these contact
  submanifolds intersect $N$ and each other transversely, i.e.\ $T_pN
  + T_pS_j = T_pM$ at every $p\in S_j\cap N$ and $T_qS_i + T_qS_j =
  T_qM$ at every $q\in S_i\cap S_j$.

  Then we find a neighborhood of $N$ in $M$ that can be represented as
  $(N\times\Disk{2},\left.\alpha\right|_{TN} + r^2\,d\phi)$ such that
  $S_j$ is in this neighborhood of the form $(S_j\cap N)\times
  \Disk{2}$.
\end{propo}
\begin{proof}
  Start by choosing a metric on $S_1\cap \dotsm \cap S_r$, and extend
  this metric first over all $S_1\cap \dotsm\cap \widehat{S_j}\cap
  \dotsm \cap S_r$, then to all $S_1\cap\dotsm\cap \widehat{S_i}\cap
  \dotsm \cap \widehat{S_j}\cap \dotsm \cap S_r$, and so on until the
  metric is defined on all $S_1,\dotsc,S_r$.  Now define finally the
  metric on the rest of $M$.  By considering the exponential map along
  $N$, we obtain a tubular neighborhood $N$ diffeomorphic to
  $N\times\Disk{2}$ such that every $S_j$ is given by $(S_j\cap
  N)\times\Disk{2}$.

  The standard neighborhood theorem (see for example
  \cite{Geiges_Handbook}) guarantees that we have a contactomorphism
  from $(S_1\cap \dotsm \cap S_r\cap N)\times \Disk{2}$ to itself,
  which deforms the contact form to
  $\left.\alpha\right|_{T(S_1\cap\dotsm \cap S_r\cap N)} +
  r^2\,d\phi$.  This map is generated by a vector field which vanishes
  on $(S_1\cap\dotsm\cap S_r\cap N)\times\{0\}$, and hence we can
  easily extend the contactomorphism to a diffeomorphism $\Phi$ from
  any $(S_1\cap \dotsm\cap \widehat{S_j}\cap \dotsm \cap S_r\cap
  N)\times \Disk{2}$ to itself such that $\Phi$ leaves $(S_1\cap
  \dotsm\cap \widehat{S_j}\cap \dotsm \cap S_r\cap N)\times\{0\}$
  fixed.  By suitably extending the vector field successively over all
  of the submanifolds and then to the rest of the manifold, we finally
  obtain a diffeomorphism that converts the contact form on $(S_1\cap
  \dotsm \cap S_r\cap N)\times \Disk{2}$ into the desired type, that
  leaves $N\times\{0\}$ pointwise fixed and all other submanifolds
  $(S_{i_1}\cap \dotsm \cap S_{i_k} \cap N) \times\Disk{2}$ invariant
  as subsets.

  Now use again the neighborhood theorem, this times for $(S_2\cap
  \dotsm \cap S_r\cap N)\times\{0\}$ in $(S_2\cap \dotsm \cap S_r\cap
  N)\times \Disk{2}$.  Note that the Moser vector field vanishes on
  $(S_1\cap \dotsm \cap S_r\cap N)$, because there the contact form
  was already brought into the desired shape in the previous step.
  Extend the vector successively over all submanifolds of
  $N\times\Disk{2}$ like before, but take care to choose it to vanish,
  if possible.  The flow of this vector field maps $N\times\Disk{2}$
  to itself, keeps $N\times\{0\}$ pointwise fixed, and all other
  submanifolds $(S_{i_1}\cap \dotsm \cap S_{i_k} \cap N)
  \times\Disk{2}$ invariant as subsets.  This map converts the contact
  form into the desired one on $(S_2\cap \dotsm \cap S_r\cap N)\times
  \Disk{2}$.

  Now repeat the process for $(S_1\cap S_3\cap \dotsm \cap S_r\cap
  N)\times \Disk{2}$.  Here care has to be taken for the Moser field
  not to destroy the form on $(S_2\cap \dotsm \cap S_r\cap N)\times
  \Disk{2}$, which was already arranged correctly in the previous
  step, but if the vector field is extended to vanish wherever
  possible, this submanifold is not moved at all by the flow.  By
  continuing with this process, one finally proves the statement in
  the proposition.
\end{proof}

Now consider a contact submanifold $S\hookrightarrow M$ that has
transverse intersection with $N$.  The neighborhood of $N$ in $M$ can
be represented as $(N\times\Disk{2},\left.\alpha\right|_{TN} +
r^2\,d\phi)$ in such a way that $S$ is of the form $(N\cap S)\times
\Disk{2}$ in this neighborhood.  Then it follows that the Presas
gluing $M \cup_N N\times\T^2$ contains the gluing $S \cup_{N\cap S}
(N\cap S)\times \T^2$ as a contact submanifold.  This construction
works for several submanifolds $S_1,\dotsc,S_r$ that satisfy the
assumptions in Proposition~\ref{umgebungssatz fuer zwei mgktn}.

\section{The proof of Corollary~\ref{hauptfolgerung}}\label{sec:hauptbeweis}

The general construction to prove Corollary~\ref{hauptfolgerung} is
considerably more complicated than the $5$--dimensional one.  The
proof works by induction.  We start with a contact sphere that
contains a $PS$--overtwisted contact submanifold of codimension~$2k$.
In each induction step, we raise the dimension of the submanifold by
two until finally the sphere itself is $PS$--overtwisted.

Let $S:=(\S^{2n-1},\alpha)$ be a contact sphere that contains
codimension~$2$ contact submanifolds $S_1,\dotsc,S_k$ (with $k\le
n-2$) with the following properties (see also Figure~\ref{bild:
  induktion}):
\begin{enumerate}
\item{} Every $S_j$ is a sphere that is unknotted in~$S$.
\item{} Every two spheres $S_i$ and $S_j$ ($i\ne j$) intersect
  transversely, and the intersection
  \begin{equation*}
    S_{i_1,\dotsc,i_r} := S_{i_1}\cap \dotsm \cap S_{i_r}
  \end{equation*}
  of any combination of these spheres is a contact $(2n - 2r
  -1)$--sphere, which is unknotted in any of the spheres
  $S_{i_1,\dotsc,\widehat i_s,\dotsc,i_r}$.
\item{} Finally the lowest dimensional sphere $S_{1,\dotsc,k}$ is
  $PS$--overtwisted.
\end{enumerate}

\begin{figure}[htbp]
  \begin{picture}(0,0)%
    \includegraphics{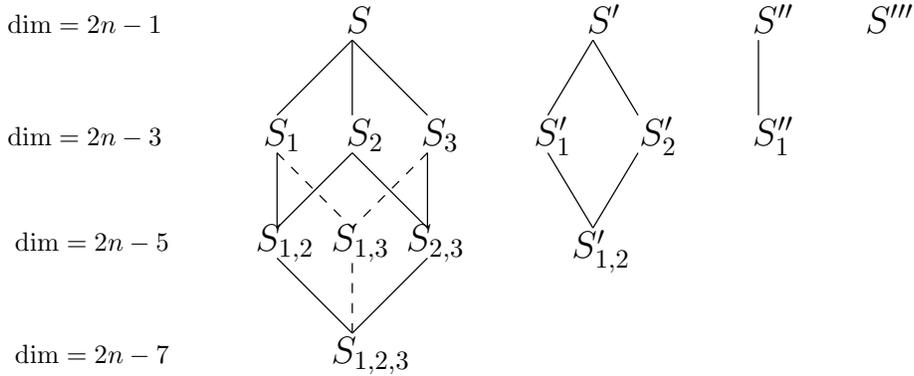}%
  \end{picture}%
  \setlength{\unitlength}{4144sp}
  \begingroup\makeatletter\ifx\SetFigFont\undefined
  \gdef\SetFigFont#1#2#3#4#5{%
    \reset@font\fontsize{#1}{#2pt}%
    \fontfamily{#3}\fontseries{#4}\fontshape{#5}
    \selectfont}%
  \fi\endgroup%
  \begin{picture}(5130,2223)(1261,-4309)

    \put(3646,-3571){\makebox(0,0)[lb]{\smash{{\SetFigFont{14}{16.8}{\familydefault}{\mddefault}{\updefault}{\color[rgb]{0,0,0}$S_{2,3}$}%
          }}}}

    \put(4726,-2266){\makebox(0,0)[lb]{\smash{{\SetFigFont{14}{16.8}{\familydefault}{\mddefault}{\updefault}{\color[rgb]{0,0,0}$S^\prime$}%
          }}}}

    \put(4411,-2941){\makebox(0,0)[lb]{\smash{{\SetFigFont{14}{16.8}{\familydefault}{\mddefault}{\updefault}{\color[rgb]{0,0,0}$S_1^\prime$}%
          }}}}

    \put(5041,-2941){\makebox(0,0)[lb]{\smash{{\SetFigFont{14}{16.8}{\familydefault}{\mddefault}{\updefault}{\color[rgb]{0,0,0}$S_2^\prime$}%
          }}}}

    \put(5716,-2941){\makebox(0,0)[lb]{\smash{{\SetFigFont{14}{16.8}{\familydefault}{\mddefault}{\updefault}{\color[rgb]{0,0,0}$S_1^{\prime\prime}$}%
          }}}}

    \put(5716,-2266){\makebox(0,0)[lb]{\smash{{\SetFigFont{14}{16.8}{\familydefault}{\mddefault}{\updefault}{\color[rgb]{0,0,0}$S^{\prime\prime}$}%
          }}}}

    \put(4636,-3616){\makebox(0,0)[lb]{\smash{{\SetFigFont{14}{16.8}{\familydefault}{\mddefault}{\updefault}{\color[rgb]{0,0,0}$S_{1,2}^\prime$}%
          }}}}

    \put(6391,-2266){\makebox(0,0)[lb]{\smash{{\SetFigFont{14}{16.8}{\familydefault}{\mddefault}{\updefault}{\color[rgb]{0,0,0}$S^{\prime\prime\prime}$}%
          }}}}

    \put(2791,-2941){\makebox(0,0)[lb]{\smash{{\SetFigFont{14}{16.8}{\familydefault}{\mddefault}{\updefault}{\color[rgb]{0,0,0}$S_1$}%
          }}}}

    \put(3286,-2941){\makebox(0,0)[lb]{\smash{{\SetFigFont{14}{16.8}{\familydefault}{\mddefault}{\updefault}{\color[rgb]{0,0,0}$S_2$}%
          }}}}

    \put(3736,-2941){\makebox(0,0)[lb]{\smash{{\SetFigFont{14}{16.8}{\familydefault}{\mddefault}{\updefault}{\color[rgb]{0,0,0}$S_3$}%
          }}}}

    \put(3286,-2266){\makebox(0,0)[lb]{\smash{{\SetFigFont{14}{16.8}{\familydefault}{\mddefault}{\updefault}{\color[rgb]{0,0,0}$S$}%
          }}}}

    \put(3196,-4246){\makebox(0,0)[lb]{\smash{{\SetFigFont{14}{16.8}{\familydefault}{\mddefault}{\updefault}{\color[rgb]{0,0,0}$S_{1,2,3}$}%
          }}}}

    \put(2746,-3571){\makebox(0,0)[lb]{\smash{{\SetFigFont{14}{16.8}{\familydefault}{\mddefault}{\updefault}{\color[rgb]{0,0,0}$S_{1,2}$}%
          }}}}

    \put(3196,-3571){\makebox(0,0)[lb]{\smash{{\SetFigFont{14}{16.8}{\familydefault}{\mddefault}{\updefault}{\color[rgb]{0,0,0}$S_{1,3}$}%
          }}}}
    \put(1261,-2266){\makebox(0,0)[lb]{\smash{{\SetFigFont{10}{12.0}{\familydefault}{\mddefault}{\updefault}{\color[rgb]{0,0,0}$\dim
              = 2n-1$}%
          }}}}
    \put(1261,-2941){\makebox(0,0)[lb]{\smash{{\SetFigFont{10}{12.0}{\familydefault}{\mddefault}{\updefault}{\color[rgb]{0,0,0}$\dim
              = 2n-3$}%
          }}}}
    \put(1306,-3571){\makebox(0,0)[lb]{\smash{{\SetFigFont{10}{12.0}{\familydefault}{\mddefault}{\updefault}{\color[rgb]{0,0,0}$\dim
              = 2n-5$}%
          }}}}
    \put(1306,-4246){\makebox(0,0)[lb]{\smash{{\SetFigFont{10}{12.0}{\familydefault}{\mddefault}{\updefault}{\color[rgb]{0,0,0}$\dim
              = 2n-7$}%
          }}}}
  \end{picture}%
  \caption{The intersection of all codimension~$2$ spheres is
    $PS$--overtwisted.  In each induction step  we increase the dimension of the $PS$--overtwisted
    submanifold by two.}\label{bild: induktion}
\end{figure}

We are then able to construct a new contact sphere $S^\prime =
(\S^{2n-1},\alpha^\prime)$ which satisfies the conditions in the above list for
$k-1$ instead of $k$.  More explicitly we mean that $S^\prime$
contains $k-1$ unknotted codimension~$2$ spheres
$S_1^\prime,\dotsc,S_{k-1}^\prime$, which are contact submanifolds,
all possible intersections are contact spheres, unknotted in any of
the other higher dimensional spheres, and $S^\prime_{1,\dotsc,k-1}$ is
$PS$--overtwisted.  By induction, we can continue these steps until we
find a contact structure on $\S^{2n-1}$, such that $\S^{2n-1}$
contains a $PS$--overtwisted unknotted $(2n-3)$--sphere.  In the next
step, we finally obtain then the $PS$--overtwisted contact structure
on the $(2n-1)$--sphere itself (the proof in Section~\ref{sec: beweis
  dimension 5} amounts to this last step).

\subsection{Start of induction}

To start the process in arbitrary dimension, consider the manifold
$S:=(\S^{2n-1},\alpha_-)\subset \C^n$ defined in
Section~\ref{sec:link-haendige sphaeren}.  Every sphere $S_j :=
\{(z_1,\dotsc,z_n)\in\S^{2n-1}|\, z_j = 0\}$ for
$j\in\{1,\dotsc,n-2\}$ is contact, and all of the possible
intersections also are.  The $3$--sphere $S_{1,\dotsc,n-2} =
\{(0,\dotsc,0,z_1,z_2)\in\S^{2n-1}\}$ is overtwisted, see Remark~\ref{rem:overtwistedforms}. Hence it is
possible to start the construction in any odd dimension~$2n-1$, and so
assume that the induction step is true for some $k<n$.  Then we have
to show that we find a $(2n-1)$--dimensional contact sphere such that
the statements are also true for~$k-1$.

\subsection{Construction of $S^\prime$ and $S^\prime_1,\dotsc,S^\prime_{k-1}$}

Apply the Presas gluing on $S = \S^{2n-1}$ along $S_k$, i.e.\
construct the manifold
\begin{equation*}
  M_0 := \S^{2n-1} \cup_{S_k} S_k\times\T^2\;.
\end{equation*}
By Section~\ref{sec: presas gluing auf untermannigfaltigkeiten}, it
contains all of the gluings on the other spheres $S_j$ with
$j=1,\dotsc,k-1$, along the intersection $S_{j,k} = S_j \cap S_k$
\begin{equation*}
  M_0^{j} := S_j \cup_{S_{j,k}} S_{j,k}\times\T^2\;,
\end{equation*}
and of course also the gluings of any other subspace
$S_{j_1,\dotsc,j_r}$ along the intersection $S_{j_1,\dotsc,j_r,k} =
S_{j_1,\dotsc,j_r}\cap S_k$.
\begin{equation*}
  M_0^{j_1,\dotsc,j_r} := S_{j_1,\dotsc,j_r} \cup_{S_{j_1,\dotsc,j_r,k}} S_{j_1,\dotsc,j_r,k}\times\T^2\;.
\end{equation*}
In particular, it follows from Theorem~\ref{satz: presas gluing} that
$M_0^{1,\dotsc,k-1}$ is $PS$--overtwisted, because we obtained it by
gluing along $S_{1,\dotsc,k}$, which is by our assumptions
$PS$--overtwisted.  By using surgery on all of these submanifolds, we
will be able to convert every $M_0^{(\dotsm)}$ into a sphere, which will prove property $(3)$ of the list at the beginning of Section~\ref{sec:hauptbeweis}.

Using that each sphere is unknotted in the next higher dimensional
one, we obtain with arguments that are completely analogous to those of
Section~\ref{sec: beweis dimension 5} that
$\pi_1(M_0^{j_1,\dotsc,j_r}) \cong \Z *\Z$, and that the homology is
$H_0(M_0^{j_1,\dotsc,j_r}) \cong H_{2n-1-2r}(M_0^{j_1,\dotsc,j_r})
\cong \Z$, $H_1(M_0^{j_1,\dotsc,j_r}) \cong
H_{2n-2-2r}(M_0^{j_1,\dotsc,j_r}) \cong \Z^2$, and all other homology
groups are trivial.

\begin{figure}[htbp]
  \begin{picture}(0,0)%
    \includegraphics{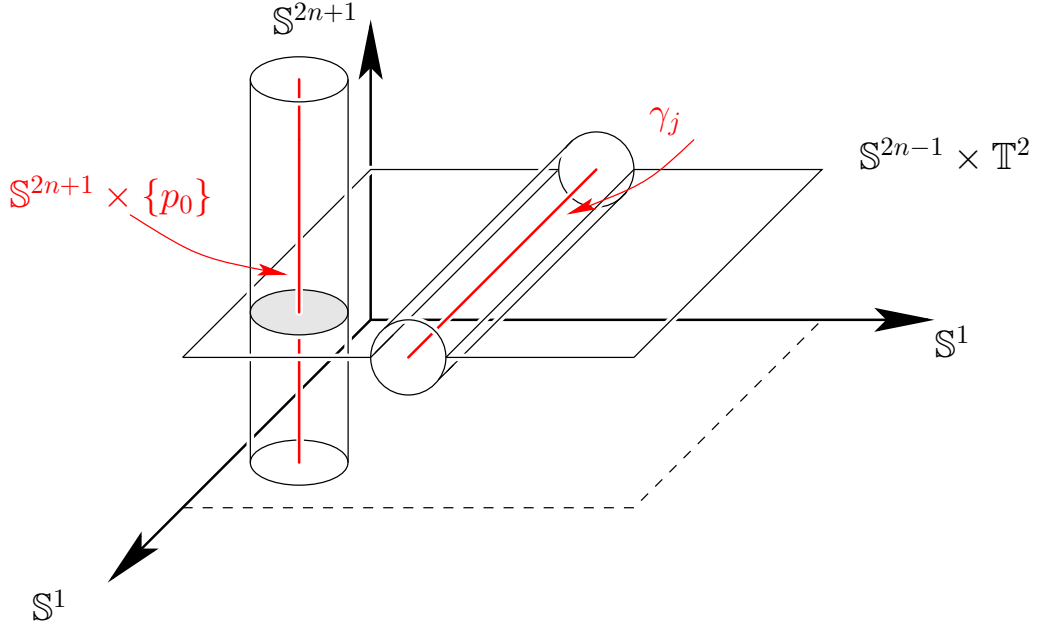}%
  \end{picture}%
  \setlength{\unitlength}{4144sp}%
  \begingroup\makeatletter\ifx\SetFigFont\undefined%
  \gdef\SetFigFont#1#2#3#4#5{%
    \reset@font\fontsize{#1}{#2pt}%
    \fontfamily{#3}\fontseries{#4}\fontshape{#5}%
    \selectfont}%
  \fi\endgroup%
  \begin{picture}(5557,3715)(1441,-5510)
    \put(1576,-5461){\makebox(0,0)[lb]{\smash{{\SetFigFont{14}{16.8}{\familydefault}{\mddefault}{\updefault}{\color[rgb]{0,0,0}$\S^1$}%
          }}}}
    \put(6976,-3886){\makebox(0,0)[lb]{\smash{{\SetFigFont{14}{16.8}{\familydefault}{\mddefault}{\updefault}{\color[rgb]{0,0,0}$\S^1$}%
          }}}}
    \put(5266,-2491){\makebox(0,0)[lb]{\smash{{\SetFigFont{14}{16.8}{\familydefault}{\mddefault}{\updefault}{\color[rgb]{1,0,0}$\gamma_j$}%
          }}}}
    \put(3016,-1951){\makebox(0,0)[lb]{\smash{{\SetFigFont{14}{16.8}{\familydefault}{\mddefault}{\updefault}{\color[rgb]{0,0,0}$\S^{2n+1}$}%
          }}}}
    \put(6526,-2761){\makebox(0,0)[lb]{\smash{{\SetFigFont{14}{16.8}{\familydefault}{\mddefault}{\updefault}{\color[rgb]{0,0,0}$\S^{2n-1}\times\T^2$}%
          }}}}
    \put(1441,-2986){\makebox(0,0)[lb]{\smash{{\SetFigFont{14}{16.8}{\familydefault}{\mddefault}{\updefault}{\color[rgb]{1,0,0}$\S^{2n+1}\times \{p_0\}$}%
          }}}}
  \end{picture}%
  \caption{Apply surgery on the curves $\gamma_1$ and $\gamma_2$.
    These two loops generate the fundamental group of all submanifolds
    considered.}\label{bild: erzeuger der fundamentalgruppe}
\end{figure}

Represent the generators of the fundamental group of the lowest
dimensional manifold $M_0^{1,\dotsc,k-1}$ by the two loops $\gamma_1,
\gamma_2$ of the form
\begin{align}
\label{eq:curves_to_contract}
  \begin{aligned}
    \gamma_1:\,&[0,1] \to S_{1,\dotsc,k}\times\T^2, \, t\mapsto
    (p_1;e^{2\pi it},1) \\
    \gamma_2:\,&[0,1] \to S_{1,\dotsc,k}\times\T^2, \, t\mapsto
    (p_2;1,e^{2\pi it})\;,
  \end{aligned} 
\end{align}
with $p_1\ne p_2\in S_{1,\dotsc,k}$ fixed.  We cannot only assume that
$\gamma_1,\gamma_2$ do not intersect the fiber
$S_{1,\dotsc,k}\times\{p_0\} $ along which we perform the Presas
gluing (Figure~\ref{bild: erzeuger der fundamentalgruppe}), but by
choosing the two points $p_1,p_2$ suitably, it can also be achieved
that the plastikstufe created in $M_0^{1,\dotsc,k-1}$ by
Theorem~\ref{satz: presas gluing} is not touched by any of the two
paths.  The conformal symplectic normal bundles of these loops are trivial, and hence
the $h$--principle guarantees that a perturbation of
$\gamma_1,\gamma_2$ in $M_0^{1,\dotsc,k-1}$ turns them into isotropic
curves, so that they can be used to apply contact surgery.  Note that
the loops $\gamma_1, \gamma_2$ chosen do not only generate
$\pi_1(M_0^{1,\dotsc,k-1})$ but also the fundamental group of any
other of the manifolds $M_0^{j_1,\dotsc,j_r}$ including the one of the
maximal space~$M_0$.

From Proposition~\ref{umgebungssatz fuer zwei mgktn}, we can easily
deduce the following corollary.

\begin{coro}
\label{cor:normalform_submanifolds}
  The neighborhood of $M_0^{1,\dots,k-1}$ is contactomorphic to
  \begin{equation*}
    \Bigl(M_0^{1,\dots,k-1} \times \C^{k-1},\,
    \left.\alpha\right|_{T M_0^{1,\dots,k-1}} +
    \frac{i}{2}\,\sum_{j=1}^{k-1}(z_j\,d\bar z_j  - \bar z_j\, d z_j)\Bigr) \;,
  \end{equation*}
  where $M_0^{j}$ is represented by $M_0^{1,\dots,k-1} \times
  \{(z_1,\dotsc,z_{k-1})\in\C^{k-1}|\, z_j = 0\}$.
\end{coro}
\begin{proof}
  Consider $M_0^1, \dotsc, M_0^{k-1}$ in $M_0$.  They satisfy the
  conditions of Proposition~\ref{umgebungssatz fuer zwei mgktn}, and
  so there is a tubular contact neighborhood of $M_0^1$
  \begin{equation*}
    \Bigl(M_0^{1} \times \C = \{(p;z_1)\},
    \left.\alpha\right|_{TM_0^1} +
    \frac{i}{2}\,(z_1\,d\bar z_1  - \bar z_1\, d z_1)\Bigr)
  \end{equation*}
  such that $M_0^j$ is given by $M_0^{1,j}\times\C$.  Repeat the step
  for $M_0^{1,2},\dotsc,M_0^{1,k-1}$ in $M_0^1$.  We obtain a
  neighborhood of the form
  \begin{equation*}
    \Bigl(M_0^{1,2} \times \C^2 = \{(p^\prime;z_1,z_2)\},
    \left.\alpha\right|_{TM_0^{1,2}} +
    \frac{i}{2} \sum_{j=1}^2(z_j\,d\bar z_j  - \bar z_j\, d z_j)\Bigr)\;.
  \end{equation*}
  This step can be iterated until one arrives at the neighborhood
  described in the corollary we want to prove.
\end{proof}

With this neighborhood theorem, we will be able to apply contact
$1$--surgery on~$M_0$ along the curves $\gamma_1$ and $\gamma_2$.

\subsubsection{Surgery compatible with submanifolds}

To describe the contact surgery, we briefly recall Weinstein's picture
for surgery \cite{WeinsteinSurgery}. Consider $\R^{2n+4}$ with
coordinates $(\vec x,\vec y,z_1,z_2,w_1,w_2)$ and symplectic form
\begin{equation*}
  \omega = d\vec x \wedge d \vec y +
  \sum_{i=1}^2 dz_i\wedge dw_i \;. 
\end{equation*}
The vector field
\begin{equation*}
  \frac{1}{2}\,(\vec x \frac{\partial}{\partial \vec x}+\vec y \frac{\partial}{\partial \vec y})+\sum_{i=1}^2 \left( 2z_i\frac{\partial}{\partial z_i} -w_i\frac{\partial}{\partial w_i} \right)
\end{equation*}
is Liouville, and it is transverse to the non-zero level sets of the
function
\begin{equation*}
  f = \frac{1}{4}\,(\vec x^2+\vec y^2)+
  \sum_{i=1}^2 \left( z_i^2-\frac{1}{2}w_i^2 \right) \;.
\end{equation*}
In particular, $f^{-1}(-1) \cong \R^{2n+2}\times \S^1$ is a contact
hypersurface, and the circle
\begin{equation*}
  \gamma_{\mathrm{model}} =\{ (0,0,0,0,w_1,w_2) ~|~w_1^2+w_2^2= 2 \}
\end{equation*}
is isotropic.  By a neighborhood theorem, the loops $\gamma_1$ and
$\gamma_2$ considered above have a neighborhood that is
contactomorphic to a neighborhood of $\gamma_{\mathrm{model}}$ in
$f^{-1}(-1)$.

\begin{figure}[htbp]
  \begin{picture}(0,0)%
    \includegraphics{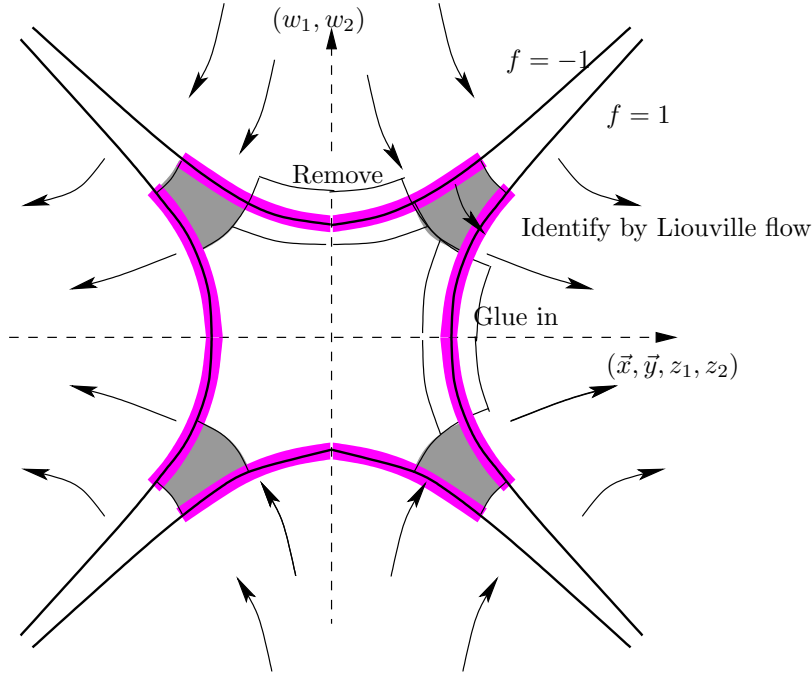}%
  \end{picture}%
  \setlength{\unitlength}{3947sp}%
  \begingroup\makeatletter\ifx\SetFigFont\undefined%
  \gdef\SetFigFont#1#2#3#4#5{%
    \reset@font\fontsize{#1}{#2pt}%
    \fontfamily{#3}\fontseries{#4}\fontshape{#5}%
    \selectfont}%
  \fi\endgroup%
  \begin{picture}(4845,4224)(364,-3373)
    \put(4126,-1486){\makebox(0,0)[lb]{\smash{{\SetFigFont{10}{12.0}{\familydefault}{\mddefault}{\updefault}{\color[rgb]{0,0,0}$(\vec x,\vec y,z_1,z_2)$}%
          }}}}
    \put(2026,689){\makebox(0,0)[lb]{\smash{{\SetFigFont{10}{12.0}{\familydefault}{\mddefault}{\updefault}{\color[rgb]{0,0,0}$(w_1,w_2)$}%
          }}}}
    \put(3594,-622){\makebox(0,0)[lb]{\smash{{\SetFigFont{10}{12.0}{\familydefault}{\mddefault}{\updefault}{\color[rgb]{0,0,0}Identify by Liouville flow}%
          }}}}
    \put(4126, 89){\makebox(0,0)[lb]{\smash{{\SetFigFont{10}{12.0}{\familydefault}{\mddefault}{\updefault}{\color[rgb]{0,0,0}$f=1$}%
          }}}}
    \put(3826,-1186){\makebox(0,0)[rb]{\smash{{\SetFigFont{10}{12.0}{\familydefault}{\mddefault}{\updefault}{\color[rgb]{0,0,0}Glue in}%
          }}}}
    \put(3496,434){\makebox(0,0)[lb]{\smash{{\SetFigFont{10}{12.0}{\familydefault}{\mddefault}{\updefault}{\color[rgb]{0,0,0}$f=-1$}%
          }}}}
    \put(2153,-294){\makebox(0,0)[lb]{\smash{{\SetFigFont{10}{12.0}{\familydefault}{\mddefault}{\updefault}{\color[rgb]{0,0,0}Remove}%
          }}}}
  \end{picture}%
  \caption{Weinstein's picture for contact
    surgery.}\label{fig_weinstein_surgery}
\end{figure}

The next step consists of gluing in the $1$--handle lying between
$f^{-1}(-1)$ and $f^{-1}(1)$.  A piece of the set $f^{-1}(-1)$ can be
identified via the Liouville flow with a contact manifold that is
close to $f^{-1}(1)$ (see Figure~\ref{fig_weinstein_surgery}).  The
precise construction that is necessary can be found in
\cite{Geiges_IntroductionBook}.  Hence we can replace a neighborhood
of the curve $\gamma_{\mathrm{model}}$ by the surgery $\Disk{2} \times
\S^{2n+1}$. As the Liouville field is transverse to the set
$f^{-1}(1)$, we see that the surgered manifold is contact.

We still need to show that we can choose the contact surgery on $M_0$
in such a way that it induces contact surgery on all the submanifolds
$M_0^J$, where $J$ is an index set. To see that the surgery can be
made compatible, we need to specify the framing more precisely, which
can be done by using an induction.
Corollary~\ref{cor:normalform_submanifolds} is of relevance here.

Let us denote the curve where we perform surgery by $\gamma$. We start
with the contact submanifold $M_0^{1,\ldots,k}$. This manifold also
contains the curve $\gamma$, which is still isotropic. Therefore we
can identify a tubular neighborhood of $\gamma$ with $f_J^{-1}(-1)$ as
in the above model. We have added the subscript $J$ to indicate the
different dimensions, i.e.\ the number of $\vec x$ coordinates depends
on the size of index set $J$. At this stage there is still some
freedom in choosing the framing. For clarity, we will also give
$\gamma$ a subindex to indicate in what submanifold we consider the
curve, i.e.~$\gamma_J$ denotes the restriction of the curve $\gamma$
to the submanifold~$M_0^J$.

Next, suppose we have fixed the framing on $M_0^{J,j}$ such that
contact surgery on $\gamma_J$ induces the desired contact surgery on
submanifolds $M_0^{J,j,K}$ in $M_0^{J,j}$ indexed by $K$. Let us now
look at $M_0^{J,j} \subset M_0^{J}$. A neighborhood of $M_0^{J,j}$ in
$M_0^J$ looks like
\begin{equation*}
  M_0^{J,j} \times \Disk{2}\;,
\end{equation*}
so we may write for the contact form by
Corollary~\ref{cor:normalform_submanifolds}
\begin{equation*}
  \alpha_{J}=\alpha_{J,j}+x\,dy-y\,dx \;,
\end{equation*}
if we use coordinates $(x,y)$ for the disk $\Disk{2}$. Note that
$\partial_x$ and $\partial_y$ trivialize the normal bundle of
$M_0^{J,j}$ in $M_0^{J}$. We use this to extend the framing of
$\gamma_{J,j}$ to a framing of $\gamma_J$ in $M_0^{J}$; we simply add
the vector fields $\partial_x$ and $\partial_y$. Note that $\partial_x$
and $\partial_y$ can be chosen to lie in $TM_0^{1,\ldots,\hat
  \j,\ldots,k}$.

To use Weinstein's picture for surgery, we use the same identifications as before, but we add a coordinate to the vectors $\vec x$ and $\vec y$. To be more precise, we identify a piece of $f_{J,j}^{-1}(-1)$ with a neighborhood of $\gamma_{J,j}$ by using the map
\begin{align*}
  \begin{aligned}
    \iota_{J,j}:~ f_{J,j}^{-1}(-1) & \to & N_{J,j}(\gamma_{J,j}) \;.
  \end{aligned}
\end{align*}

Note that
\begin{equation*}
  f_{J}(\vec x,0,\vec y,0,z_1,z_2,w_1,w_2)=f_{J,j}(\vec x,\vec y,z_1,z_2,w_1,w_2) \;.
\end{equation*}
We can extend the contactomorphism $\iota_{J,j}$ to a contactomorphism
$\iota_{J}$ such that
\begin{equation*}
  \iota_{J}(\vec x,0,\vec y,0,z_1,z_2,w_1,w_2)=
  \iota_{J,j}(\vec x,\vec y,z_1,z_2,w_1,w_2) \;.
\end{equation*}
This shows that we can choose the surgery compatible with the
submanifold structure, i.e.~surgery on $\gamma \subset M_0$ induces
the desired surgery on the submanifolds $M_0^J$ for all index sets
$J\subset \{ 1,\dotsc, k\}$.

Let us now return to our previous notation where we have chosen curves
$\gamma_1$ and $\gamma_2$ in $M_0$ as given in
Formula~\eqref{eq:curves_to_contract}.  Denote the surgered manifolds
$M_0$ and $M_0^1,\dotsc,M_0^{k-1}$ by $S^\prime$ resp.\ by
$S^\prime_1,\dotsc, S^\prime_{k-1}$.  By what we said above, the
intersection
\begin{equation*}
  S^\prime_{j_1,\dotsc,j_r} := S_{j_1} \cap \dotsm \cap S_{j_r}
\end{equation*}
is equal to $M_0^{j_1,\dotsc,j_r}$ surgered along the curves
$\gamma_1$ and $\gamma_2$. The manifolds $S'_J$ are homeomorphic to
spheres by the Poincaré conjecture as proved by Smale. Furthermore
they intersect transversely, so this almost proves property~$(2)$ of
the list. There is, however, one missing part, namely we still need to
show that the spheres are in fact diffeomorphic to the standard
sphere. This shall be done in the next section.

\subsection{Every manifold $S^\prime_{j_1,\dotsc,j_r}$ is diffeomorphic
  to the standard sphere}

Our next aim is to see that all the manifolds $M_0^{j_1,\dotsc,j_r}$
are converted by the surgery into smooth spheres.  For this, we will
show that each of these surgered manifolds is the boundary of a ball.

Consider the manifold $H$ consisting of the thickened $2$--torus
$\Disk{2n}\times \T^2$ where we attach a $2$--handle along
$\{p_1\}\times \S^1\times \{1\}$, another one along $\{p_2\}\times
\{1\} \times \S^1 $ and a $(2n)$--handle along
$(\partial\Disk{2n})\times \{p_0\}$ (with $p_1 \ne p_2\in
  \partial\Disk{2n}$, $p_0\in \T^2 - (\S^1\times\{1\} \cup \{1\}\times
  \S^1)$).  If we have chosen above the right framing for attaching
  the handles, then the boundary of $H$ is diffeomorphic to a manifold
  $S_{j_1,\dotsc,j_r}^\prime$ which was defined as
  $M_0^{j_1,\dotsc,j_r}$ surgered along $\gamma_1$ and $\gamma_2$. So $H$ provides a topological filling for $S_{j_1,\dotsc,j_r}^\prime$. 
  Using a further surgery on the interior of $H$, we will convert $H$
  into a ball.  This will show that every $S_{j_1,\dotsc,j_r}^\prime$ is
  diffeomorphic to the standard sphere.

The fundamental group of the handle body $H$ constructed above is
trivial, because $\pi_1(\Disk{2n}\times \T^2) \cong \Z^2$, attaching
the $(2n)$--handle does not change the fundamental group, and the
final two $2$--handles then kill $\pi_1(H)$.  The homology of $H$ can
be computed with a Mayer-Vietoris sequence, where we set $A=
\Disk{2n}\times \T^2$ and $B= \Disk{2n} \dot\cup \Disk{2} \dot\cup
\Disk{2}$ and $A\cap B = \S^{2n-1}\times\{p_0\} \dot\cup
\{p_1\}\times\S^1\times\{1\} \dot\cup \{p_2\}\times\{1\}\times \S^1$,
and so
\begin{equation*}
  0 \to H_2(A) \oplus H_2(B)\to H_2(H) \to
  H_1(A\cap B) \to H_1(A) \oplus H_1(B) \to 0 \;,
\end{equation*}
which simplifies to
\begin{equation*}
  0 \to \Z \to H_2(H) \to \Z^2 \to \Z^2 \to 0 \;.
\end{equation*}
This means that $H_2(H)\cong \Z$.  Hence the manifold $H$ is obviously not
a $(2n+2)$--ball.  For higher homology the sequence gives
\begin{equation*}
  0 \to H_{2n}(H) \to H_{2n-1}(A\cap B) \to 0 \;,
\end{equation*}
so that $H_{2n}(H) \cong \Z$.  All other homology groups $H_k(H)$ with
$2<k<2n$ are trivial.  The generator of $H_{2n}(H)$ can be represented
by the cycle composed of the $(2n)$--handle of $H$ and the
$(2n)$--disk $\Disk{2n}\times\{p_0\}$.
  
To convert $H$ into the desired $(2n+2)$--ball, note that the
generator of $H_2(H)$ can be represented by an embedded $2$--sphere.
On an abstract level this is clear by the Hurewicz theorem, but for
later observations we will need to find the generator explicitly.  The
$2$--torus $\{0\}\times\T^2$ at the core of the thickened torus with
which the construction of $H$ began, clearly generates $H_2(H)$ as can
be seen from the Mayer-Vietoris sequence.  Of course this torus can be
moved to $\{p_1\}\times \T^2$.  One of the $2$--handles has been
attached to $\{p_1\}\times\S^1\times\{1\}$.  We will cut out the
annulus
\begin{equation*}
  A = \bigl\{\{p_1\}\times\S^1\times\{e^{it}\}\bigm|\,
  t\in(-\epsilon,\epsilon)\bigr\} \;,
\end{equation*}
and attach two disks to $\partial A$ which lie in the $2$--handle.
Depending on the framing we used to attach the $2$--handle, the two
boundary circles $\partial A_\pm$ of $A$ are given by
\begin{equation*}
  \S^1\to \Disk{2n}\times\partial\Disk{2},\, 
  e^{i\phi}\mapsto (\pm \z(\phi),e^{i\phi}) \;.
\end{equation*}
This curves can be easily extended to the interior of the $2$--handle
by writing
\begin{equation*}
  \Disk{2}\to \Disk{2n}\times\Disk{2},\, 
  re^{i\phi}\mapsto (\pm r\z(\phi),re^{i\phi}) \;.
\end{equation*}
By a general position argument, the two disks can also be made
disjoint, because the lowest dimension for $H$ we are considering
is~$6$.  The manifold obtained from the torus by cutting out $A$,
gluing the two disks and smoothing out the corners is an embedded
$2$--sphere, which represents the generator of $H_2(H)$, because the
annulus $A$ together with the two disks represents a boundary in the
chain complex of the $2$--handle.

To perform surgery along this sphere, we have to show that its normal
bundle in $H$ is trivial.  All rank~$N$ vector bundles over the
$2$--sphere can be constructed by taking trivializations over the
hemispheres and then gluing both parts over the equator.  The gluing
operation is described by a map $\partial\Disk{2}\to \SO(N)$, and any
two homotopic maps give rise to isomorphic vector bundles.  Hence the
vector bundles are classified by $\pi_1(\SO(N)) \cong \Z_2$ if $N>2$,
and both bundles can be distinguished by the second Stiefel-Whitney
class.  Here we see that $\iota^*w_2(TH) = w_2(\left.TH\right|_{\S^2})
= w_2(T\S^2) + w_1(T\S^2)\,w_1(\nu \S^2) + w_2(\nu\S^2) =
w_2(\nu\S^2)$.  On the $2$--sphere, we can obtain
$w_2(\left.TH\right|_{\S^2})$ by evaluating $w_2$ over the fundamental
class of $\S^2$. We compute
\begin{align*}
  \pairing{w_2(\left.TH\right|_{\S^2})}{[\S^2]} &=
  \pairing{\iota^*w_2(TH)}{[\S^2]} =
  \pairing{w_2(TH)}{\iota_*[\S^2]} \\
  &= \pairing{w_2(TH)}{[\{0\}\times\T^2]} = 0 \;,
\end{align*}
where we have used that the evaluation map is independent of the
representative used for the homology class, and that $TH$ can be
trivialized over the $2$--torus.

Perform now surgery on this $2$--sphere, i.e.~cut out $\S^2 \times
\Disk{2n}$ and glue in $\Disk{3}\times \S^{2n-1}$ and denote the
resulting manifold by $\widetilde H$.  The manifold constructed this
way is simply connected (A closed loop $\gamma$ in $\widetilde H$ can
be made disjoint from the surgery region such that it represents a
loop in $H$, which can be contracted without intersecting the surgery
region.)

The Mayer-Vietoris sequence gives with the notation $B = \S^2 \times
\Disk{2n}$, $A= H - B$, $A\cup B = H$, and $A\cap B = \S^2\times
\S^{2n-1}$ that $H_k(A) = \{0\}$ for $2<k<2n-1$, and with
\begin{equation*}
  0 \to \Z \to H_2(A)\oplus \Z \to \Z \to 0 \to H_1(A)
  \to 0 \;,
\end{equation*}
one sees that $H_2(A)\cong \Z$.  Below it will be important to
understand the map $\iota_*:\, H_2(A\cap B) \to H_2(A)$.  In the
sequence above, the generator of $H_2(A\cap B)$ goes to $(a,1)\in
H_2(A)\oplus H_2(B)$ with $a\in \Z$.  The map $H_2(B) \to H_2(H)$
sends generator to generator, so that $a$ has to be a generator,
because $(a,1)$ lies in the kernel of the map $\iota_A - \iota_B$ in
the sequence above.  It follows that $\iota_*:\, H_2(A\cap B) \to
H_2(A)$ is an isomorphism.

For the higher groups compute the rest of the Mayer-Vietoris sequence:
\begin{equation*}
  0 \to \Z \to H_{2n+1}(A) \to 0
  \to 0 \to H_{2n}(A) \to \Z \to \Z \to H_{2n-1}(A) \to 0 \;.
\end{equation*}
It follows that $H_{2n+1}(A) \cong \Z$, $H_{2n}(A)$ is either trivial
or isomorphic to $\Z$, and $H_{2n-1}(A)$ can be trivial or isomorphic
to either $\Z_p$ or $\Z$.  To recognize that these two groups are
trivial, we will study the map $H_{2n}(H) \to H_{2n-1} (\S^2\times
\S^{2n-1})$ and prove that it is an isomorphism.  The generator of
$H_{2n}(H)$ can be written as the core of the $(2n)$--handle of $H$
composed with the ball $\Disk{2n}\times \{p_0\} \subset
\Disk{2n}\times \T^2$.  To see the image under the connecting
homomorphism we need to represent the generator as the sum of two
chains, one which lies in $A$ and one which lies in $B$.  The boundary
of these chains lies in $A\cap B$ (because they cancel each other),
and the boundary of one of these chains gives a representative for the
image of the generator under the connecting homomorphism.  In our
situation $B$ is a neighborhood of the $2$--sphere which generates
$H_2(H)$, and as we showed above, this $2$--sphere is obtained by
taking the torus $\{0\}\times \T^2 \subset \Disk{2n}\times \T^2$ and
attaching two disks, which lie in one of the $2$--handles.  The set
$B$ corresponds thus to a tubular neighborhood of $\S^2$, which
coincides outside a small neighborhood of the $2$--handle with a
tubular neighborhood of $\{0\}\times\T^2$.  The intersection of this
tubular neighborhood with the generator of $H_{2n}(H)$ gives a small
disk lying in a fiber of the normal bundle of $\S^2$. The boundary of
the small disk represents the generator of $H_{2n-1}(B)$ as we wanted
to show.  Hence the connecting homomorphism is a bijection and both
groups $H_{2n}(A)$ and $H_{2n-1}(A)$ are trivial.

Use now the notation $A = H - B$, $\widetilde B= \Disk{3} \times
\S^{2n-1}$, $A\cap \widetilde B = \S^2 \times \S^{2n-1}$, and
$\widetilde H = A\cup \widetilde B$.  One sees immediately that
$H_k(\widetilde H) = \{0\}$ for all $2n-1 > k > 3$, and
\begin{equation*}
  \dotsm \to H_3(\widetilde H) \to H_2(A\cap \widetilde B) \to H_2(A)
  \to H_2(\widetilde H)  \to 0 \;.
\end{equation*}
We showed that the map $H_2(A\cap B)\cong\Z \to H_2(A)\cong\Z$ is an
isomorphism, hence it follows that $H_2(\widetilde H) \cong \{0\}$.
The higher parts of the Mayer-Vietoris sequence give
\begin{align*}
  0 \to \Z \to H_{2n+1}(A) \to H_{2n+1}(\widetilde H)
  \to & 0 \\
  0 \to H_{2n}(A) \to H_{2n}(\widetilde H) \to \Z \to H_{2n-1}(A)
  \oplus \Z \to H_{2n-1}(\widetilde H) \to & 0 \;.
\end{align*}
In the first sequence, we can use that the first map is an
isomorphism, as could be seen from the Mayer-Vietoris sequence when
doing the first step of the surgery, and hence $H_{2n+1}(\widetilde
H)$ is trivial.  Since we know that $H_{2n}(A)$ and $H_{2n-1}(A)$ are
trivial, the second sequence simplifies to
\begin{equation*}
  0 \to H_{2n}(\widetilde H) \to \Z
  \to \Z \to H_{2n-1}(\widetilde H) \to 0 \;.
\end{equation*}
The map in the middle is also an isomorphism, and we obtain that
$H_{2n}(\widetilde H) \cong H_{2n-1}(\widetilde H) \cong \{0\}$.
To compute $H_3(\widetilde H)$ analyze the sequence
\begin{equation*}
  \dotsm \to H_4(\widetilde H) \to H_3(A\cap \widetilde B) \to H_3(A)\oplus H_3(\widetilde B)
  \to H_3(\widetilde H)  \to 0 \;.
\end{equation*}
The case $n=2$ is different from the case $n>2$, but $H_3(\widetilde
H)$ vanishes for both.

We just have shown that the homology of $\widetilde H$ is equal to the
one of a point, and it is also simply connected, hence $\widetilde H$
is diffeomorphic to the ball $\Disk{2n+2}$ (\cite{MilnorHCobordism}),
and so the boundary $\partial\widetilde H \cong
S_{j_1,\dotsc,j_r}^\prime$ is a smooth standard sphere.

Together with our previous remarks, we have now established properties $(2)$ and $(3)$ from the list for the manifolds $S_J'$. In last section, we will prove the remaining property of the list to conclude the induction step.

\subsection{The sphere $S_{j_1,\dotsc,j_r,j}'$  is unknotted in $S_{j_1,\dotsc,j_r}'$ }
Throughout this section, we shall use the index set $J=\{j_1,\dotsc,j_r\}$ to abbreviate the notation. In order to show that the spheres are unknotted, we shall use the following characterization of unknots due to Levine,
\cite{Levine_knots}
\begin{theorem}[Levine]
  A smoothly embedded sphere $\iota:\,\S^{k-2}\to \S^k$ is unknotted
  if and only if $\pi_1(\S^k-\iota(\S^{k-2}))\cong \Z$.
\end{theorem}

This theorem reduces the problem to a computation of
the fundamental group of the complement of $S_{J,j}'$ in
$S_{J}'$. It is helpful to first consider the fundamental
group before the $1$--surgery. In that case we have
\begin{equation*}
  M_0^{J}=S_J \cup_{S_{J,k}}S_{J,k}\times \T^2 \;,
\end{equation*}
and we have a similar expression for $M_0^{J,j}$.  Next, we note that
$S_{J,k}$ and $S_{J,j,k}$ are, by induction, unknotted, so we can
write
\begin{equation*}
  S_J = \S^1\times \Disk{2n}\cup \Disk{2}\times \S^{2n-1} 
  \text{ and }
  S_{J,j} = \S^1\times \Disk{2n-2}\cup \Disk{2}\times \S^{2n-3} \;. 
\end{equation*}
These
decompositions are such that $S_{J,k}=\{ 0 \}\times \S^{2n-1}\subset
S_J$ and $S_{J,j,k}=\{ 0 \}\times \S^{2n-3}\subset S_{J,j}$. Hence we
have decompositions that are adapted to the fiber connected sum. In
other words, we can write
\begin{equation*}
  M_0^{J} = \S^1\times \Disk{2n} \cup_{A\times S_{J,k}} S_{J,k}\times
  (\T^2-\{p\}) \;,
\end{equation*}
where we glue along the annulus times fiber $A\times S_{J,k}$. We
proceed by computing the fundamental group of the complement of
$M_0^{J,j}$ in $M_0^{J}$.  We have
\begin{equation*}
  M_0^{J} - M_0^{J,j} = \S^1\times (\Disk{2n}-\Disk{2n-2})
  \cup_{A\times (S_{J,k}-S_{J,j,k})} (S_{J,k}-S_{J,j,k})\times
  (\T^2-\{p\}) \;.
\end{equation*}
Here we glue along an annulus times fiber ($S_{J,k}-S_{J,j,k}$ in this
case).  The Seifert-Van~Kampen theorem gives us the fundamental group
\begin{align*}
  \pi_1( M_0^{J}-M_0^{J,j} ) &= \langle a,b,c,d,e~|~ab=ba,~cd=dc,~ce=ec,~ded^{-1}e^{-1}=a,~b=c \rangle \\
  &\cong \langle c,d,e|~cd=dc,~ce=ec \rangle\cong \Z\oplus \left( \Z * \Z \right)
\end{align*}
Here the generators $a,b$ denote the two commuting generators of the
left-hand side, $\S^1\times (\Disk{2n}-\Disk{2n-2})$.  The generator
$c$ can be represented by a curve in $S_{J,k}-S_{J,j,k}$ that
generates $\pi_1(S_{J,k}-S_{J,j,k})\cong \Z$. Finally, the elements
$d$ and $e$ can be represented by a longitudinal curve and a meridian
of $\T^2-\{p\}$.  Because of the product structure of the right-hand
side, the two generators $d,e$ commute with $c$.  The other two
relations come from the amalgamation in the Seifert-Van~Kampen
theorem. We see in particular that the fundamental group of the
complement of $M_0^{J,j}$ has already a simple structure before the
surgery.

Now take any element $W$ in $\pi_1(S_{J}'-S_{J,j}' ,q)$, where $q$
denotes the basepoint.  By a homotopy, we can arrange $W$ to be
represented by a product of curves that lie outside the surgery
region, i.e.\ we assume that the curves lie in $M_0^{J}-M_0^{J,j}$.
When thought of as an element in $\pi_1(M_0^{J}-M_0^{J,j},q)$ our
above computation shows that such a product of curves can be written
as a word in $c,d$ and $e$.  Note however that a curve $\Gamma$
representing $d$ or $e$ can be written in terms of $c$ after the
surgery.  Indeed, we can homotope such a curve $\Gamma$ close to the
curves $\gamma_1$ or $\gamma_2$ from
Equation~\eqref{eq:curves_to_contract}.  More precisely, for dimension
reasons we can find a homotopy $H$ that is disjoint from the surgery
locus such that
\begin{align*}
  H(0,t) &= \Gamma(t) \\
  H(1,t) &= \tilde \gamma_1(t)=(\tilde p;e^{2\pi i t},1)\in \S^n \times
  \T^2
\end{align*}
in case $\Gamma$ is homotopic to $\gamma_1$ as a curve in $M_0^{J}$.
We have a similar homotopy $H$ for the case that $\Gamma$ is homotopic
to $\gamma_2$ as a curve in $M_0^{J}$.  The homotopy $H$ takes place in
$M_0^{J}-M_0^{J,j}$ and only homotopes $\Gamma$ to a curve $\tilde
\gamma_i$ that is close to $\gamma_i$ in the following sense.

We have chosen $\tilde p$ disjoint from a tubular neighborhood
$N_\epsilon$ of $S_{1,\ldots,k}$ that we use to perform surgery, but
inside the larger tubular neighborhood $N_{2\epsilon}$ such that
$\tilde p$ is close to either $p_1$ or $p_2$ as given by
Equation~\eqref{eq:curves_to_contract}.  Now we see that $\Gamma$ can
be simplified. Follow the homotopy $H$ and then push the curve $\tilde
\gamma_i$ into the surgered region.  Before the surgery, this region
(with $M_0^{J,j,k}$ removed) looks like
\begin{equation}
  \label{eq:surgery_region}
  (\Disk{n}-\Disk{n-2})\times \S^1 \;.
\end{equation}
When we homotope $\tilde \gamma_i$ into this region, $\tilde \gamma_i$
can arranged to have the form
\begin{equation*}
  t\mapsto (f(t);e^{2\pi i t}) \;.
\end{equation*}
The precise form of $f$ depends on the chosen framing of the
neighborhood of $\S^1$.  Now the surgery replaces the
set~\eqref{eq:surgery_region} by
\begin{equation*}
  (\S^{n-1}-\S^{n-3})\times \Disk{2} \;.
\end{equation*}
Hence we can homotope $\tilde \gamma_i$ to a curve of the form
\begin{equation*}
  (f(t);1) \;,
\end{equation*}
which represents $c^k$ for some $k$. This establishes our claim.

We see therefore that $\pi_1(S_{J}'-S_{J,j}' ,q)$ can be presented by
a group with at most one generator.  Since we know that
$H_1(S_{J}'-S_{J,j}')\cong \Z$, we have
\begin{equation*}
  \pi_1(S_{J}'-S_{J,j}',q)\cong \Z \;.
\end{equation*}
As a result of Levine's criterion we obtain that the embedding
$S_{J,j}'$ into $S_{J}'$ is unknotted.  This proves property $(1)$ of
the list and finishes the induction step.

\bibliographystyle{amsalpha}

\bibliography{main}


\end{document}